\documentclass[10pt]{article}
\usepackage{amsmath}
\usepackage{amssymb}
\usepackage{amsfonts}
\numberwithin{equation}{section}
\usepackage[dvips]{graphicx}

\newcommand{\qed}{\vrule width 5pt height 5pt depth 0pt}
\newtheorem{thm}{THEOREM}[section]

\newtheorem{lemma}[thm]{LEMMA}

\def\today{{\ifcase\month\or
 January\or February\or March\or April\or May\or June\or
 July\or August\or September\or October\or November\or December\fi
 \space\number\day, \number\year}}

\newcommand{\calD}{{\mathcal D}}

\def\Re{\operatorname{Re}}    
              
\def\C{{\mathbb C}}       
           \def\calE{\mathcal E}

\def\R{{\mathbb R}}

               \def\bS{\mathbf s}
\def\bX{\mathbf x}       \def\bW{\mathbf w}
\def\Vhat{\widehat V}   \def\Vtilde{\widetilde V}   

\def\cont{{\calD(\varrho,\theta)}}    
\def\decis{{\Gamma(\varrho,\theta)}}   
     
\newcommand{\domain}{{\Pi_+}}

\newcommand{\heat}{\boxplus}

\def\lap{{\mathcal L}}   
\def\Vhat{\widehat V}

\newcommand{\diff}{{\mathrm d}}

\def\sgn{\operatorname{sgn}}

\begin{document}

\font\titlefont=cmr10 scaled \magstep3
\centerline{\titlefont{On the asymptotic free boundary for}}
\centerline{\titlefont{the American put option problem}}

\bigskip

\centerline{\sl H\aa kan Hedenmalm}

\centerline{\sl at the Royal Institute of Technology, Stockholm}

\bigskip\bigskip

\noindent{\bf Abstract.} In practical work with American put options, it is
important to be able to know when to exercise the option, and when not to
do so. In computer simulation based on the standard theory of geometric 
Brownian motion for simulating stock price 
movements, this problem is fairly easy to handle for options with a short 
lifespan, by analyzing binomial trees. It is considerably more challenging 
to make the decision for American put options with long life\-span. In 
order to provide a satisfactory analysis, we look at the corresponding free 
boundary problem, and show that the free boundary -- which is the curve that 
separates the two decisions, to exercise or not to -- has an asymptotic 
expansion, where the coefficient of the main term is expressed as an 
integral in terms of the free boundary. This raises the perspective that 
one could use numerical simulation to approximate the integral and thus 
get an effective way to make correct decisions for long life options.   

\section{Introduction}\label{sec1}

\subsection{General background}

\noindent{\bf Initial remarks.}
The standard model of stock price movement, as proposed by Samuelson, 
is geometric Brownian motion. On the basis of this model, it is possible 
to analyze prices of derivative securities such as options, by use of the 
arbitrage principle, which postulates that a riskless portfolio must, in 
the absence of transaction costs, earn the risk-free interest rate. In this 
paper we shall be concerned with the American put option, which is a 
contract that allows the holder to sell a stock at a fixed price -- 
independent of market movements -- at any moment during the duration of 
the contract. It should be pointed out that the academic community is not 
unilaterally in favor of the geometric Brownian motion model (see, for 
instance, \cite{Mand}).
\medskip

\noindent{\bf The payoff function; geometric Brownian motion.}
Let $t$ be a time parameter, {\sl which expresses the time remaining until
the option expires}. We note that $t$ flows backwards with time, so that $t$ 
decreases as real time passes by. Typically, we are
interested in an interval $0\le t\le T$, where $t=T$ corresponds to
the time when the option is issued and $t=0$ is the time of
expiration. {\sl It is convenient to express all money in terms of its  
equivalent value at the deadline $t=0$}. Let $\bS(t)$ denote the stock price 
at time $t$. If $r$ is the risk-free continuously compounded interest rate, 
assumed to be constant, then $\bS(t)\exp(-rt)$ is the nominal value (that is, 
the dollar amount we would see on the screen) of the stock at (remaining) 
time $t$. For reasons of convenience, {\sl we normalize the nominal exercise 
price of the put option to equal $1$}. 
When we think of the stock price as a free parameter, we write 
$s$ in place of $\bS(t)$. The {\sl reward function} is 
\begin{equation}
V_r(t,s)=\max\big\{0,\exp(rt)-s\big\};
\label{eq-rewardf}
\end{equation}
it expresses the payoff earned by exercising the option contract at the 
point $(t,s)$. We consider put options of the {\sl American}
type, which means that the holder is at liberty to exercise at any time from 
purchase at $t=T$ to expiration $t=0$. Let $\Vhat_r(t,s)$ denote the correct 
price of the option at (remaining) time $t$. It is well known from the 
Black-Scholes analysis that $\Vhat_r(t,s)$ is obtained by optimizing -- over 
all stopping strategies -- the expected value of the reward function over 
all stock paths starting at $\bS(t)=s$, under the assumption of risk 
neutrality. Risk neutrality means that the expected growth of a risky asset 
like $\bS(t)$ is postulated to equal that of a riskless one (and since all 
monetary values are discounted to their equivalents at time $t=0$, the 
expected growth is $0$). The assumption of geometric Brownian motion leads 
to the infinitesimal equation
$$\diff\bS(t)=\lambda\,\bS(t)\,\diff t+\sigma\, \bS(t)\,\diff\bW(-t),$$
where $\sigma^2$ is the variance per unit of time, $\lambda$ a drift rate (the
intrinsic growth rate of the stock), and $\bW(-t)$ is the unit Brownian 
motion (we write $-t$ to indicate that our time parameter flows backwards). 
The postulate of risk neutrality translates into $\lambda=0$. 
By Ito's formula, this then leads to
\begin{equation}
\diff\big(\log \bS(t)\big)=\frac{\sigma^2}2\,\diff t+\sigma\, \diff\bW(-t).
\label{1-1}
\end{equation}
\medskip

\subsection{The obstacle problem}

\noindent{\bf Shift to ordinary Brownian motion.} It is more convenient to 
work with ordinary Brownian motion rather than geometric Brownian motion, 
and so we introduce the stochastic process
$$\bX(t)=\frac{\sigma}{\sqrt{2}}\,t-\frac{\sqrt{2}}{\sigma}\,\log \bS(t),$$
which follows a standard Wiener process, as can be seen from (\ref{1-1}):
$$\diff\bX(t)=-\sqrt{2}\,\diff\bW(-t).$$
Likewise, we switch from the pair $(t,s)$ to $(t,x)$ as our basic 
coordinates, where
$$x=\frac{\sigma}{\sqrt{2}}\,t-\frac{\sqrt{2}}{\sigma}\,\log s$$
is the free variable corresponding to $\bX(t)$.
\rm 
In the new coordinate system, the reward function (\ref{eq-rewardf}) takes 
the form
\begin{equation}
V_{r,\sigma}(t,x)=\max\left\{0,e^{rt}-e^{\sigma^2t/2-\sigma x/\sqrt{2}}
\right\}.
\label{1-2}
\end{equation}
Let
$$\heat=\frac{\partial}{\partial t}-\frac{\partial^2}{\partial
x^2}$$
be the heat operator. We say that a function $f$ of the variables
$(t,x)$ is {\sl caloric} if $\heat\,f=0$, {\sl subcaloric} if it is 
real-valued and $\heat\,f\le0$, and {\sl supercaloric} if it is 
real-valued and $\heat\,f\ge0$.
We need to introduce the upper half-plane (or positive time plane)
$$\R_+^2=\big\{(t,x):\,t\in\R_+,\,x\in\R\big\},$$
where $\R$ is the real line, and $\R_+=\{t\in\R:\,t>0\}$.
The following statement is rather well-known, and may be taken as the formal 
definition of the envelope function $\Vhat$.  

\begin{thm} The envelope function $\Vhat_{r,\sigma}$ is supercaloric in 
$\R_+^2$. Moreover, it equals the infimum of all $C^2$-smooth supercaloric 
functions that majorize $V_{r,\sigma}(t,x)$ on $\R_+^2$.
\label{def-thm-1}
\end{thm}


The scaling properties of the heat operator $\heat$ permit us to reduce the
number of parameters by $1$. We introduce the real parameter $\varrho$, 
$-1<\varrho<1$, and let $V_\varrho$ denote $V_{r,\sigma}$, with the parameter 
settings 
$$r=1-\varrho^2,\quad\sigma=\sqrt{2}(1+\varrho):$$
\begin{equation}
V_\varrho(t,x)=\max\Big\{0,e^{(1-\varrho^2)t}-
e^{(1+\varrho)^2 t-(1+\varrho)x}\Big\}.
\label{1-3}
\end{equation}
Likewise, $\Vhat_\varrho$ denotes the function $\Vhat_{r,\sigma}$ with the 
same settings. If, however, $\sigma$ and $r$ are given, we pick
$\alpha,\varrho$ according to
\begin{equation*}
\alpha=\frac{\sigma^2+2r}{2\sqrt{2}\,\sigma},\qquad
\varrho=\frac{\sigma^2-2r}{\sigma^2+2r},
\end{equation*}
and recover the function $V_{r,\sigma}$ from the formula
\begin{equation}
V_{r,\sigma}(t,x)=V_\varrho(\alpha^2t,\alpha x);
\label{eq-recover-obst}
\end{equation}
by the scaling properties of $\heat$, we then recover the envelope as well:
\begin{equation}
\Vhat_{r,\sigma}(t,x)=\Vhat_\varrho(\alpha^2t,\alpha x).
\label{eq-recover-env}
\end{equation}
In the sequel, we shall only be concerned with the function $V_\varrho$; 
moreover, we shall drop the subscript $\varrho$ whenever this does not lead 
to confusion.
\medskip

\noindent{\bf An affine shift of coordinates.} If a function $f(t,x)$ is 
caloric, {\sl then so is the transformed function}
\begin{equation}
e^{\beta^2t+\beta x}\,f(t,x+2\beta t);
\label{subst-1}
\end{equation}
here, $\beta$ is a real parameter. The calculation that shows this {\sl also
reveals that the supercaloric functions are preserved under the
transformation}. Actually, it is possible to give a complete characterization
of the transformations of this type that preserve the caloric functions.

Choose $\beta=\varrho$ in the substitution (\ref{subst-1}), and introduce
the function
$$V'_\varrho(t,x)=e^{\varrho^2t+\varrho x}\,V_\varrho(t,x+2\varrho t).$$
Then the least supercaloric majorant to $V'_\varrho$, denoted by 
$\Vhat'_\varrho$, is related to $\Vhat_\varrho$ in a straightforward fashion:
$$\Vhat'_\varrho(t,x)=e^{\varrho^2t+\varrho x}\,
\Vhat_\varrho(t,x+2\varrho t).$$
In other words, {\sl we may as well replace $V_\varrho$ by $V'_\varrho$ in 
our considerations. }
The function $V'_\varrho$ is simpler-looking:
$$V'_\varrho(t,x)=\max\left\{0,e^{t+\varrho x}-e^{t-x}\right\},
\qquad (t,x)\in\R^2_+.$$
\sl In the sequel, we shall consider only the transformed function 
$V'_\varrho$, and write $V_\varrho$ for it. \rm 
\medskip

\noindent{\bf Introduction of a new parameter.} We introduce the parameter 
$\theta$, confined to $0<\theta<+\infty$, and let 
$V_{\varrho,\theta}$ be the function
$$V_{\varrho,\theta}(t,x)=
e^{t+\varrho x}-\frac12\big(1+\theta-\varrho+\theta\varrho\big)\,
e^{t-x}-\frac12(1-\theta)(1+\varrho)\,e^{t+x},\quad (t,x)\in\Pi_+,$$
extended to the whole positive time half-plane $\R_+^2$ by
$$V_{\varrho,\theta}(t,x)=0,\qquad (t,x)\in\R^2_+\setminus\Pi_+.$$
We readily calculate that on $\R_+^2$,
\begin{equation}
\heat\,V_{\varrho,\theta}(t,x)
=-\theta\,(1+\varrho)\,e^t\,\delta_0(x)+(1-\varrho^2)\,
e^{t+\varrho x}\,1_{\domain}(t,x),
\label{eq-heatV}
\end{equation}
where $1_E$ denotes the characteristic function of the set $E$, and 
$\delta_0$ is the unit Dirac mass at $0$. 

From this, it is evident that the role of $\theta$ is to scale the mass
distribution along the $t$-axis. The value $\theta=1$ corresponds to the
put option problem in the introduction. We have the corresponding envelope 
function $\Vhat_{\varrho,\theta}$, and Theorem \ref{def-thm-1} generalizes 
to the new setting.

\subsection{The free boundary}\label{subsec-1.3}

\noindent\bf The continuation region. \rm The region 
$$\calD(\varrho,\theta)=\Big\{(t,x)\in\R^2_+:\,V_{\varrho,\theta}(t,x)
<\Vhat_{\varrho,\theta}(t,x)\Big\}$$
is called the {\sl continuation region}, while the boundary curve
$$\Gamma(\varrho,\theta)=\partial\calD(\varrho,\theta)\cap\R^2_+$$
is the {\sl free boundary}, or {\sl decision boundary}. Note that 
$$\R^2_+\setminus \Pi_+\subset\calD(\varrho,\theta).$$
We parametrize $\Gamma(\varrho,\theta)$ by
$$x=\phi_{\varrho,\theta}(t),\qquad 0<t<+\infty,$$
and note that it is well-known that $\phi_{\varrho,\theta}(t)$ is an 
increasing function, with $\phi_{\varrho,\theta}(0)=0$ and 
\begin{equation}
\lim_{t\to+\infty}\phi_{\varrho,\theta}(t)=\mu(\varrho,\theta)
=\frac{1}{1+\varrho}\,
\log\left(1+\theta\,\frac{1+\varrho}{1-\varrho}\right),
\label{def-mu}
\end{equation}
where the equation is taken to define the constant $\mu=\mu(\varrho,\theta)$.
Moreover, at least for $\theta=1$, it is known that the function 
$x=\phi_{\varrho,\theta}(t)$ is continuous and concave 
(see, for instance, \cite{EE}).
In financial terms, in the continuation region, we keep the options contract,
while at the decision boundary, we exercise it.
The behavior of the free boundary near $t=0$ has been studied extensively
(see \cite{BBRS}):
$$\phi_{\varrho,1}(t)\sim(1+\varrho)\sqrt{2t\log\frac1t},
\qquad t\to0\quad(\text{for }\,\theta=1).$$
Here, we shall focus on the behavior of the free boundary as $t\to+\infty$.
We shall demonstrate that $\phi_{\varrho,\theta}(t)$ has the asymptotic
expansion
$$\phi_{\varrho,\theta}(t)=\mu(\varrho,\theta)-\beta_1\,t^{-3/2}\,e^{-t}
-\ldots-\beta_N\,t^{-N-1/2}\,e^{-t}+O_{L^p}\big(t^{-N-1}\,e^{-t}\big),$$
where $\beta_j=\beta_j(\varrho,\theta)$ are certain real-valued coefficients.
The term 
$$O_{L^p}\big(t^{-N-1}\,e^{-t}\big)$$ 
stands for an $L^p(\R_+)$-function times $t^{-N-1}\,e^{-t}$, and we are free
to choose $p$ with $2<p<+\infty$. 
We also show how to express $\beta_j$ in terms of an integral, which should
permit numerical algorithms to yield good approximations of $\beta_j$. In
particular, 
$$\beta_1=\frac{e^{-\varrho\mu}}{2\sqrt{\pi}}\int_0^{+\infty}
\phi_{\varrho,\theta}(t)\,
e^{\varrho\,\phi_{\varrho,\theta}(t)}\,\phi'_{\varrho,\theta}(t)\,
e^t\,\diff t,$$
which we may -- by integration by parts -- also write in the form
\begin{multline*}
\beta_1=\frac{1}{2\sqrt{\pi}}\bigg\{
\frac{\theta\mu\, e^{-\varrho\mu}}{1-\varrho}+
\frac{(1+\varrho\mu)\,e^{-\varrho\mu}-1}{\varrho^2}\\
+\frac1{\varrho^2}\int_0^{+\infty}\Big[
\big(\varrho\mu-\varrho\phi_{\varrho,\theta}(t)+1\big)
e^{\varrho\,(\phi_{\varrho,\theta}(t)-\mu)}
-1\Big]\,e^t\,\diff t\bigg\},
\end{multline*}     
since it is a consequence of our analysis that
$$\int_0^{+\infty}\phi'_{\varrho,\theta}(t)\,e^{\varrho\,
\phi_{\varrho,\theta}(t)}\,e^t\,\diff t=
\frac{\theta}{1-\varrho}.$$ 

\subsection{Asymptotic estimation of the free boundary}

\noindent{\bf Estimation of $\widehat V_{\varrho,\theta}$ from above and 
below.} 
We need to get a grasp of the asymptotic behavior of the function 
$\phi_{\varrho,\theta}(t)$ as $t\to+\infty$. 
Let $\Vtilde_{\varrho,\theta}$ be the solution to the 
heat equation with boundary values equal to $V_{\varrho,\theta}$ on the two 
half-lines
$$\big\{(t,x)\in\R^2:\,\,\,t=0,\quad -\infty<x<\mu\big\}$$
and
$$\big\{(t,x)\in\R^2:\,\,\,0<t<+\infty,\quad x=\mu\big\};$$
it is unique under mild growth restrictions. This function arises from
the stopping strategy of stopping on the two given lines, and hence we
must have 
$$\Vtilde_{\varrho,\theta}(t,x)\le\Vhat_{\varrho,\theta}(t,x),
\qquad(t,x)\in\Pi_\mu,$$ 
where
\begin{equation}
\Pi_\mu=\big\{(t,x)\in\R^2:\,\,\,0<t<+\infty,\,\,\, -\infty<x<\mu\big\},
\label{def-Pimu}
\end{equation}
and $\mu=\mu(\varrho,\theta)$ is as in (\ref{def-mu}).

Let $\Vhat_\infty$ be the function
$$\Vhat_\infty(t,x)=\eta\,e^{t+x},\qquad (t,x)\in\R_+^2,$$
where $\eta$ is the constant
$$\eta=\frac12\,(1+\varrho)\big[e^{(\varrho-1)\mu}-1+\theta\big];$$
then
$$\Vhat_{\varrho,\theta}(t,x)\le\Vhat_\infty(t,x),
\qquad(t,x)\in\Pi_\mu.$$
A calculation reveals that
\begin{equation}
\Vhat_\infty(t,x)-V_{\varrho,\theta}(t,x)=
\frac12(1-\varrho^2)\,e^{\varrho\mu}\,
e^t(x-\mu)^2+O\big(e^t(x-\mu)^3\big),
\label{diff-1}
\end{equation}
near the line $x=\mu$. Let $V_1$ denote the function
$$V_1(t,x)=\Vhat_\infty(t,x)-
\Vtilde_{\varrho,\theta}(t,x),\qquad (t,x)\in\Pi_\mu,$$
and note that it is caloric and that it vanishes along the boundary vertical
line $x=\mu$. By the reflection principle for caloric functions
\cite[pp. 115--116]{W}, then, $V_1$ extends to a caloric function throughout
the positive time half-plane, with 
$$V_1(t,x)\equiv-V_1(t,2\mu-x),\qquad (t,x)\in\R^2_+.$$
The boundary values of $V_1$ along the $x$-axis are given by
$$V_1(0,x)=\Vhat_\infty(0,x)-\Vtilde_{\varrho,\theta}(0,x)=
\Vhat_\infty(0,x)-V_{\varrho,\theta}(0,x),\qquad x\in]-\infty,\mu].$$
This is not as explicit as desired, so we decompose
$$V_1(t,x)=V_2(t,x)-V_3(t,x),$$
where $V_2$ and $V_3$ are caloric in the positive time half-plane, 
with boundary values
$$V_2(0,x)=\eta\,\sgn(\mu-x)\,\min\big\{e^x,e^{2\mu-x}\big\},$$
and
$$V_3(0,x)=V(x)=e^{\varrho x}-\frac12(1-\varrho)e^{\varrho\mu}\,e^{\mu-x}
-\frac12(1-\theta)(1+\varrho)e^{x},\qquad x\in[0,\mu],$$
with $V_3(0,x)=0$ for $x\in\R\setminus[0,2\mu]$, and $V_3(x)=-V_3(2\mu-x)$
for $x\in[\mu,2\mu]$. The first boundary moment of $V_2$ is
$$\int_{-\infty}^{+\infty}(x-\mu)\,V_2(0,x)\,\diff x=-2\eta,$$
and that of $V_3$ is
\begin{multline*}
\int_{-\infty}^{+\infty}(x-\mu)\,V_3(0,x)\,\diff x=-2\,\frac{e^{\varrho\mu}-
1-\varrho\mu}{\varrho^2}\\+
(1-\varrho)e^{(\varrho+1)\mu}\big(e^{-\mu}-1+\mu\big)
+(1-\theta)(1+\varrho)\big(e^{\mu}-1-\mu\big),
\end{multline*}
so that
\begin{multline}
\int_{-\infty}^{+\infty}(x-\mu)\,V_1(0,x)\,\diff x=-2\eta+
2\,\frac{e^{\varrho\mu}-1-\varrho\mu}{\varrho^2}\\-
(1-\varrho)e^{(\varrho+1)\mu}\big(e^{-\mu}-1+\mu\big)
-(1-\theta)(1+\varrho)\big(e^{\mu}-1-\mu\big)<0;
\label{eq-wtmn}
\end{multline}
after all, the function $V_1(0,x)$ is positive for $x<\mu$ and negative for
$x>\mu$. We need the following observation.

\begin{lemma}
Suppose $f$ is a continuous odd complex-valued function on the real line,
which decays exponentially rapidly at infinity. Then its caloric extension
$f(t,x)$ to the positive time half plane, as given by the formula
$$f(t,x)=
\frac1{2\sqrt{\pi t}}\int_{-\infty}^{+\infty}e^{-\frac{(x-\xi)^2}{4t}}
\,f(\xi)\,\diff\xi,$$
has the asymptotics
$$f(t,x)=\frac{x}{4\sqrt{\pi}t^{3/2}}\int_{-\infty}^{+\infty}f(\xi)
\,\xi\,\diff\xi+O\big(x\,t^{-5/2}\big)$$
as $t\to+\infty$ and $x$ is kept inside a compact interval of the real line.
\end{lemma}

The proof is left as an exercise to the reader.

In view of the lemma and (\ref{diff-1}), we obtain
\begin{multline*}
\Vtilde_{\varrho,\theta}(t,x)-V_{\varrho,\theta}(t,x)
=\Vhat_\infty(t,x)-V_{\varrho,\theta}(t,x)-V_1(t,x)=
\frac12(1-\varrho^2)\,e^{\varrho\mu}\,
e^t(x-\mu)^2\\
-\frac{x-\mu}{4\sqrt{\pi}\,t^{3/2}}\int_{-\infty}^{+\infty}
V_1(0,\xi)\,(\xi-\mu)\,\diff\xi+O\big(e^t(x-\mu)^3\big)+
O\big((x-\mu)\,t^{-5/2}\big),
\end{multline*}
so that $V_{\varrho,\theta}(t,x)<\Vtilde_{\varrho,\theta}(t,x)$ holds in a 
domain of the type
$$x-\mu<-B_1\,t^{-3/2}e^{-t}+O\big(t^{-5/2}e^{-t}\big),$$
where the positive constant $B_1$ is given by
\begin{equation}
B_1=\frac{e^{-\varrho\mu}}{2\sqrt{\pi}(1-\varrho^2)}\int_{-\infty}^{+\infty}
(\mu-\xi)\,V_1(0,\xi)\,\diff\xi,
\label{def-B1}
\end{equation}
which is evaluated in equation (\ref{eq-wtmn}).
We have obtained the following statement.

\begin{lemma} The function $\phi_{\varrho,\theta}$ describing the decision 
boundary enjoys the following estimate:
\begin{equation*}
\mu-B_1\,t^{-3/2}\,e^{-t}+O\big(t^{-5/2}e^{-t}\big)\le\phi_{\varrho,\theta}(t)
<\mu,\qquad 0<t+\infty,
\end{equation*}
where the positive constant $B_1$ is given by (\ref{def-B1}).
\label{lm-1}
\end{lemma}

Note that this proves the assertion that the vertical line $x=\mu$ is an
asymptote for the decision boundary.

\section{The balayage equation}

\subsection{Derivation of the balayage equation}

\noindent{\bf Integration by parts and balayage.} 
We introduce the function
$$U_{\varrho,\theta}(t,x)=\Vhat_{\varrho,\theta}(t,x)-
V_{\varrho,\theta}(t,x);$$
in the continuation region $\cont$, it solves the overdetermined problem
$$\begin{cases}
\heat\, U_{\varrho,\theta}=-\heat V_{\varrho,\theta}
\qquad\text{on}\,\,\,\,\cont,\\
U_{\varrho,\theta}=0\qquad\qquad\quad\text{on}\,\,\,\,\partial\cont,\\
\partial_x U_{\varrho,\theta}=0\qquad\qquad\text{on}\,\,\,\,\decis.
\end{cases}$$
In the complement $\R^2_+$, we have
$$\Vhat_{\varrho,\theta}(t,x)=V_{\varrho,\theta}(t,x),\qquad
(t,x)\in\R_+^2\setminus\cont,$$
which makes $U_{\varrho,\theta}=0$ there. Also, as $U_{\varrho,\theta}=0$ on 
$\partial\cont$, we can extend the function $U_{\varrho,\theta}$ continuously 
to all of $\R^2$ by declaring $V=0$ throughout $\R^2\setminus\cont$. Since 
$U_{\varrho,\theta}$ vanishes together with its gradient along $\decis$, 
it actually solves
\begin{equation}
\left\{\begin{gathered}
\!\!\!\!\!\!\heat\, U_{\varrho,\theta}
=-1_\cont\,\heat\,V_{\varrho,\theta}\qquad\text{on}
\,\,\,\,\R^2,\\
\, U_{\varrho,\theta}=0\qquad\qquad\qquad
\text{on}\,\,\,\,\R^2\setminus\cont.
\end{gathered}
\right.
\label{3-1}
\end{equation}

Written out more explicitly using (\ref{eq-heatV}), (\ref{3-1}) assumes the 
following form:
\begin{equation}
\heat\,U_{\varrho,\theta}=\left[\theta\,(1+\varrho)\,\delta_0(x)-(1-\varrho^2)
\,1_{\domain\cap\cont}(t,x)\right]\,e^{t+\varrho x},
\label{eq-heatU}
\end{equation}
on $\R^2_+$, while
\begin{equation}
U_{\varrho,\theta}(t,x)=0,\qquad(t,x)\in\R^2\setminus\cont.
\label{eq-heatU-1}
\end{equation}
By integration by parts, if $h$ is a compactly supported $C^\infty$-smooth 
function in $\R^2$, then 
\begin{equation}
\int_{\R^2}\heat\, U_{\varrho,\theta}(t,x)\,h(t,x)\,\diff t\diff x=
\int_{\R^2}U_{\varrho,\theta}(t,x)\,\heat^*\! h(t,x)\,\diff t\diff x,
\label{eq-intbyparts}
\end{equation}
where 
$$\heat^*=-\frac{\partial}{\partial t}-\frac{\partial^2}{\partial
x^2}$$
is the adjoint heat operator. We would like to plug the functions
$$h(t,x)=h_z(t,x)=e^{-tz^2+xz}$$
into (\ref{eq-intbyparts}), for $z\in\C$, because they are 
all $\heat^*$-caloric. This function, unfortunately, 
is never compactly supported; however, it might be possible to approximate 
it by compactly supported functions so that (\ref{eq-intbyparts}) holds for 
$h=h_z$ in the limit. As we proceed in this manner, taking into account 
the known growth properties of $V_{\varrho,\theta}$, we find that we should 
restrict the complex parameter $z$ to 
\begin{equation}
1<\Re z\quad\text{and}\quad 1<\Re(z^2).
\label{restr-zeta}
\end{equation}
We compute
$$\int_{\R^2_+}\theta(1+\varrho)\,\delta_0(x)\,e^{t+\varrho x}\,
e^{-tz^2+xz}\,\diff t\diff x=\frac{\theta(1+\varrho)}{z^2-1},$$
while
\begin{multline*}
\int_{\R^2_+}(1-\varrho^2)\,1_{\Pi_+\cap\calD(\varrho,\theta)(t,x)}\,
e^{t+\varrho x}\,e^{-tz^2+xz}\,\diff t\diff x=\\
\frac{1-\varrho^2}{\varrho-z}\int_0^{+\infty}
e^{(\varrho+z)\phi_{\varrho,\theta}(t)}\,
e^{(1-z^2)t}\diff t-\frac{1-\varrho^2}{(\varrho+z)(z^2-1)}.
\end{multline*}
It now follows from (\ref{eq-heatV}) and (\ref{eq-intbyparts}) that
\begin{equation}
\int_0^{+\infty}e^{\phi_{\varrho,\theta}(t)[\varrho+z]}\,e^{-[z^2-1]t}
\,\diff t=
\frac{1-\varrho+\theta\,(\varrho+z)}{(1-\varrho)(z^2-1)}.
\label{eq-BE}
\end{equation}
We shall call this the {\sl balayage equation} for the free boundary equation
$x=\phi_{\varrho,\theta}(t)$. The reason is that the positive mass 
concentrated on the half-line
$$\big\{(t,x)\in\R^2_+:\,x=0\big\}$$
should be counterbalanced by a corresponding negative mass spread evenly over
$$\domain\cap\cont=\big\{(t,x)\in\R^2_+:\,0<x<\phi_{\varrho,\theta}(t)
\big\}.$$

It is possible to show that, under reasonable restrictions, the balayage
equation characterizes the free boundary. We shall, however, not pursue this
matter here.

\subsection{Equivalent formulations of the balayage equation}

\noindent{\bf Rewriting the balayage equation.} 
In the context of the balayage equation (\ref{eq-BE}), we introduce the 
complex variable $s=z^2$. We use the principal branch of the square root
to define $\sqrt{s}$, so that if $s$ has positive real part, then so does 
$\sqrt{s}$. In particular, if $z$ meets (\ref{restr-zeta}), then
$\sqrt{s}=z$. It follows that we may rewrite (\ref{eq-BE}) in the form
\begin{equation}
\int_0^{+\infty}e^{\phi_{\varrho,\theta}(t)[\varrho+\sqrt{s}]}\,e^{-[s-1]t}
\,\diff t=
\frac{1-\varrho+\theta\,(\varrho+\sqrt{s})}{(1-\varrho)(s-1)},
\label{eq-BE-2}
\end{equation}
for all $s\in\C$ with $\Re s>1$.

We introduce the function $\varphi_{\varrho,\theta}$,
\begin{equation*}
\varphi_{\varrho,\theta}(t)=\mu-\phi_{\varrho,\theta}(t),
\qquad t\in[0,+\infty[,
\end{equation*}
where $\mu=\mu(\varrho,\theta)$ is as before. The function 
$\varphi_{\varrho,\theta}$ is decreasing, with 
$$\varphi_{\varrho,\theta}(0)=\mu$$ 
and 
\begin{equation}
\varphi_{\varrho,\theta}(t)=O\big(t^{-3/2}e^{-t}\big)\quad\text{as}\,\,\,
t\to+\infty,
\label{eq-asymp}
\end{equation}
in view of Lemma \ref{lm-1}.
In terms of this function, the balayage equation (\ref{eq-BE-2}) assumes the 
form
\begin{equation}
\int_0^{+\infty}e^{-\varphi_{\varrho,\theta}(t)\big[\varrho+\sqrt{s}\big]}\,
e^{-(s-1)t}\,\diff t
=\frac{e^{-\mu\varrho}}{1-\varrho}\,
\frac{1-\varrho+\theta\,[\varrho+\sqrt{s}]}{s-1}\,e^{-\mu\sqrt{s}},
\label{eq-BE-3}
\end{equation}
for $\Re s>1$. {\sl This is the version we shall use many times in 
the sequel.}
By the way, an integration by parts man\oe{}uvre applied to (\ref{eq-BE-2}) 
results in the simpler-looking identity
\begin{equation*}
\int_0^{+\infty}\phi'_{\varrho,\theta}(t)\,
e^{\phi_{\varrho,\theta}(t)\big[\varrho+\sqrt{s}\big]}\,
e^{-(s-1)t}\,\diff t=\frac{\theta}{1-\varrho},
\end{equation*}
valid for $\Re s>1$. In terms of the function $\varphi_{\varrho,\theta}$, 
the relationship reads
\begin{equation}
\int_0^{+\infty}
\big|\varphi'_{\varrho,\theta}(t)\big|\,e^{-\varphi_{\varrho,\theta}(t)
\big[\varrho+\sqrt{s}\big]}\,
e^{-(s-1)t}\,\diff t=\frac{\theta\,e^{-\mu\varrho}}{1-\varrho}\,
e^{-\mu\sqrt{s}},
\label{3-8.1}
\end{equation}
where, again, $\Re s>1$.

\subsection{A general scheme for analyzing the balayage equation}

\noindent{\bf Taylor's formula for the exponential function.} 
The exponential function has a Taylor expansion about the origin:
$$e^{z}=1+z+\frac{z^2}{2!}+\ldots+\frac{z^N}{N!}+E_{N+1}(z),$$
where the remainder term $E_{N+1}(z)$ can be expressed in the form
$$E_{N+1}(z)=\frac1{N!}\int_0^z (z-\xi)^N\,e^\xi\,\diff\xi=\frac{z^{N+1}}{N!}
\int_0^1 (1-t)^N\,e^{tz}\,\diff t,$$
which represents an entire function, and enjoys the estimate
\begin{equation}
\big|E_{N+1}(z)\big|\le\frac{|z|^{N+1}}{(N+1)!}\max\big\{1,e^{\Re z}\big\}.
\label{est-rem}
\end{equation}
As we apply the above formula to the balayage equation (\ref{eq-BE-3}),
the result is
\begin{multline}
\frac1{s-1}-\big[\varrho+\sqrt{s}\big]\,\lap[\varphi_{\varrho,\theta}](s-1)
+\frac1{2!}\big[\varrho+\sqrt{s}\big]^2\,
\lap[\varphi^2_{\varrho,\theta}](s-1)-\ldots\\
+\frac{(-1)^N}{N!}\big[\varrho+\sqrt{s}\big]^N\,
\lap[\varphi^N_{\varrho,\theta}](s-1)
+\int_0^{+\infty}E_{N+1}
\Big(-\big[\varrho+\sqrt{s}\big]\varphi_{\varrho,\theta}(t)\Big)
\,e^{-(s-1)t}\,\diff t\\
=\frac{e^{-\mu\varrho}}{1-\varrho}\,
\frac{1-\varrho+\theta\,[\varrho+\sqrt{s}]}{s-1}\,e^{-\mu\sqrt{s}},
\label{3-8.3}
\end{multline}
for $\Re s>1$. Here, $\lap$ denotes the Fourier-Laplace transform, as
defined by
$$\lap[f](z)=\int_{0}^{+\infty}e^{-tz}\,f(t)\,\diff t,$$
wherever the integral converges.
The intention is to use the identity (\ref{3-8.3}) to sucessively obtain 
more information regarding the function $\varphi_{\varrho,\theta}$. 
Indeed, in view of Lemma \ref{lm-1}, we have some input to initiate the 
iterative process.

\subsection{Asymptotics of the free boundary}

\noindent{\bf Analytic continuation of the Laplace transform.}
We recall the definition of the function $\varphi$,
$$\phi_{\varrho,\theta}(t)=\mu-\varphi_{\varrho,\theta}(t),
\qquad 0<t<+\infty,$$
which is a decreasing function with $\varphi_{\varrho,\theta}(0)=\mu$ and 
the asymptotic bound (\ref{eq-asymp}). 
We use the representation (\ref{3-8.3}) with $z=\sqrt{s}$ and $N=1$ in the 
form
\begin{multline}
\lap[\varphi_{\varrho,\theta}](z^2-1)=\frac{1}{z^2-1}-
\frac{e^{-\mu\varrho}}{1-\varrho}\,\,
\frac{1-\varrho+\theta(z+\varrho)}{z^2-1}\,\,e^{-\mu z}\\
-\int_0^{+\infty}E_2\Big(-[\varrho+z]\varphi_{\varrho,\theta}(t)\Big)
\,e^{-(z^2-1)t}\,\diff t.
\label{4-3}
\end{multline}
It is immediate from the definition of $\mu=\mu(\varrho,\theta)$ that
the function
$$G(z)=\frac{1}{z^2-1}-
\frac{e^{-\mu\varrho}}{1-\varrho}\,\,
\frac{1-\varrho+\theta(z+\varrho)}{z^2-1}\,\,e^{-\mu z}$$
extends analytically to $\C\setminus\{-1\}$, and that the singularity at
$z=-1$ is a simple pole.
Moreover, by Lemma \ref{lm-1}, the left hand side converges to a holomorphic 
function for $z$ with $\Re[z^2]>0$. Furthermore, by (\ref{est-rem}), the 
function
$$\calE_2(z)=
\int_0^{+\infty}E_2\Big(-[\varrho+z]\varphi_{\varrho,\theta}(t)\Big)
\,e^{-(z^2-1)t}\,\diff t$$
expresses an even analytic function in the simply connected domain 
$\Re[z^2]>-1$, which extends continuously to the closed region 
$\Re[z^2]\ge-1$. It follows that the right hand side of (\ref{4-3}) expresses 
an analytic function in $\Re[z^2]>-1$ with the exception of a simple pole at
$z=-1$.   

We introduce the notation
$$\Phi(z)=\lap[\varphi_{\varrho,\theta}](z^2-1)$$
for $z\in\C$ with $\Re z>0$ and $\Re [z^2]>0$, and let $\Phi$ denote
the possible analytic continuation of this function beyond its initial
domain of definition. Then (\ref{4-3}) reads
\begin{multline}
\Phi(z)=\frac{1}{z^2-1}-
\frac{e^{-\mu\varrho}}{1-\varrho}\,\,
\frac{1-\varrho+\theta(z+\varrho)}{z^2-1}\,\,e^{-\mu z}\\
-\int_0^{+\infty}E_2\Big(-[\varrho+z]\varphi_{\varrho,\theta}(t)\Big)
\,e^{-(z^2-1)t}\,\diff t,
\label{4-3-1}
\end{multline}
for all $z\in\C$ with $\Re[z^2]>-1$ and $z\neq-1$, and $\Phi(z)$ is 
holomorphic in $\Re[z^2]>-1$ with the exception of a simple pole at $z=-1$,
and continuous up to the boundary.
\medskip

\noindent{\bf Estimate of the Laplace transform.} 
We need some size control of $\Phi(z)$, in order to invoke
inverse Laplace transformation and obtain information about 
$\varphi_{\varrho,\theta}$. By integration by parts,
\begin{multline}
\frac1{z+\varrho}\int_0^{+\infty}
E_2\Big(-[\varrho+z]\varphi_{\varrho,\theta}(t)\Big)\,
e^{-[z^2-1]t}\,\diff t
=\frac1{z^2-1}\left[\mu-\frac{1-e^{-\mu(z+\varrho)}}{z+\varrho}\right]\\
-\frac1{z^2-1}
\int_0^{+\infty}\varphi'_{\varrho,\theta}(t)\,
E_1\Big(-\varphi_{\varrho,\theta}(t)[z+\varrho]\Big)
\,e^{-[z^2-1]t}\,\diff t.
\label{4-4.1}
\end{multline}
Combining (\ref{4-3}) and (\ref{4-4.1}), we arrive at 
\begin{multline}
\Phi(z)
=\frac1{z^2-1}\,\left[\mu-\frac{\theta\,e^{-\mu(z+\varrho)}}{1-\varrho}
\right]\\
-\frac1{z^2-1}\int_0^{+\infty}\varphi'_{\varrho,\theta}(t)\,
E_1\Big(-\varphi_{\varrho,\theta}(t)[z+\varrho]\Big)\,e^{-[z^2-1]t}\,\diff t.
\label{4-4.2}
\end{multline}
Let us, as a matter of convenience, introduce also the function
$$\Psi(z)=\mu-[z^2-1]\,\Phi(z),$$
which is a Fourier-Laplace transform as well:
$$\Psi(z)=\lap[-\varphi'_{\varrho,\theta}](z^2-1),$$
for $z\in\C$ with $\Re z>0$ and $\Re[z^2]>0$. In terms of $\Psi(z)$, 
(\ref{4-4.2}) reads
\begin{equation}
\Psi(z)=
\frac{\theta\,e^{-\mu(z+\varrho)}}{1-\varrho}+
\int_0^{+\infty}\varphi'(t)\,
E_1\Big(-\varphi_{\varrho,\theta}(t)[z+\varrho]\Big)
\,e^{-[z^2-1]t}\,\diff t,
\label{4-4.3}
\end{equation}
for $z\in\C$ with $\Re[z^2]>-1$. 
By letting $z\to0$ along the positive real axis, we obtain from (\ref{4-4.3})
that
\begin{equation*}
\Psi(0)=
\int_0^{+\infty}|\varphi'_{\varrho,\theta}(t)|\,e^{t}\,\diff t
=\frac{\theta\,e^{-\mu\varrho}}{1-\varrho}
-\int_0^{+\infty}|\varphi'_{\varrho,\theta}(t)|\,
E_1\big(-\varphi_{\varrho,\theta}(t)\varrho\big)
\,e^{t}\,\diff t,
\end{equation*}
which simplifies to
\begin{equation}
\int_0^{+\infty}|\varphi'_{\varrho,\theta}(t)|\,
e^{-\varphi_{\varrho,\theta}(t)\varrho}\,e^{t}\,\diff t=
\frac{\theta\,e^{-\mu\varrho}}{1-\varrho}.
\label{4-4.4.1}
\end{equation}
From this and the fact that the Laplace transform of a positive function 
is dominated by its behavior along the real line, we get, by integration by
parts, that
\begin{equation}
\big|\lap[\varphi_{\varrho,\theta}](s)\big|=
O\bigg(\frac1{|s|}\bigg)\qquad\hbox{as }\,\,\,
|s|\to+\infty,\,\,\,\Re s>-1.
\label{eq-boundlap}
\end{equation}
This leads to an estimate of $\Phi(z)$ in the quarter-plane $\Re z>0$, 
$\Re[z^2]>0$.
We need size control of $\Phi(z)$ in the bigger region $\Re[z^2]>-1$. 
From (\ref{4-4.3}) and the estimate (\ref{est-rem}), we obtain
\begin{multline*}
|\Psi(z)|\le
\frac{\theta\,e^{-\mu(\varrho+\Re z)}}{1-\varrho}\\+
|z+\varrho|\int_0^{+\infty}|\varphi'_{\varrho,\theta}(t)|\,
\varphi_{\varrho,\theta}(t)\,
\max\big\{1,e^{-\varphi_{\varrho,\theta}(t)[\varrho+\Re z]}\big\}
\,e^{2t}\,\diff t.
\end{multline*}
In view of the derived integrability properties of $\varphi$ and its 
derivative, this gives the following growth estimate for $\Phi$:
\begin{equation}
|\Phi(z)|=O\left(\frac1{|z|}\,
\max\big\{1,e^{-\mu\Re z}\big\}\right),\quad\hbox{as}\,\,\,
|z|\to+\infty,\,\,\,\Re[z^2]>-1.
\label{4-4.7}
\end{equation}
\medskip

\noindent{\bf The Laplace transform and a slit domain.} 
It follows from the derived properties of $\Phi(z)$ that the Laplace transform
$\lap[\varphi_{\varrho,\theta}](s)$ extends analytically to the region 
$\Re s>-2$ minus the 
slit $]-2,-1]$ on the real line, and it extends continuously to all boundary
points with $\Re s=-2$, except for $s=-2$.  
At $s=-1$, $\lap[\varphi_{\varrho,\theta}](s)$ has a square root branch 
point, and for real $x$, $-2\le x\le-1$, we can speak of the functions 
$\lap[\varphi_{\varrho,\theta}](x+i0)$ and
$\lap[\varphi_{\varrho,\theta}](x-i0)$ as continuous limits from above and 
below. This makes the following expression well-defined:
\begin{equation*}
\Lambda(x)=\frac{i}{2\pi}\Big(\lap[\varphi_{\varrho,\theta}](-x+i0)-
\lap[\varphi_{\varrho,\theta}](-x-i0)\Big),\qquad 1\le x\le2.
\end{equation*}
Using the identity (\ref{4-4.2}), we can write this function as
\begin{multline}
\Lambda(x)=\frac{1}{\pi x}\bigg\{
\frac{\theta\,e^{-\mu\varrho}}{1-\varrho}
\,\sin\Big[\mu\sqrt{x-1}\Big]\\
-\int_0^{+\infty}|\varphi'_{\varrho,\theta}(\tau)|\,e^{-\varphi
(\tau)\varrho}\sin\Big[\varphi_{\varrho,\theta}(\tau)\sqrt{x-1}\Big]\,
e^{x\tau}\,\diff\tau
\bigg\},\qquad 1\le x\le2,
\label{deflambda}
\end{multline}
from which we see that it is real-valued on the interval in question, and
has a square root type singularity at the point $x=1$. We now extend the 
function $\Lambda$ to the whole real line by setting it equal to $0$ off
the interval $[1,2]$. We shall need the function $\varphi_\Lambda$, the 
Laplace transform of $\Lambda$:
\begin{equation*}
\varphi_\Lambda(t)=\lap[\Lambda](t)=\int_1^{+\infty}\Lambda(x)\,e^{-tx}\,
\diff x.
\end{equation*}
It is real-valued, and, due to the noted property 
$$\Lambda(x)=O(\sqrt{x-1})\qquad\text{as}\,\,\,x\to1^+,$$ 
it has the following asymptotics:
\begin{equation*}
\varphi_\Lambda(t)=O\big(t^{-3/2}e^{-t}\big),\qquad\text{as}\,\,\, 
t\to+\infty.
\end{equation*}  
The Laplace transform of $\varphi_\Lambda$ equals
\begin{equation*}
\lap[\varphi_\Lambda](s)=\frac{1}{2\pi i}\int_{-2}^{-1}\frac{
\lap[\varphi_{\varrho,\theta}](x+i0)-\lap[\varphi_{\varrho,\theta}](x-i0)}
{x-s}\,\diff x,
\end{equation*}
which means that the difference
\begin{equation*}
\lap[\varphi_{\varrho,\theta}](s)-\lap[\varphi_\Lambda](s)
\end{equation*}
is holomorphic throughout $\Re s>-2$, with continuous boundary values,
except possibly for a logarithmic singularity at $s=-2$. It has the
decay rate
\begin{equation*}
\big|\lap[\varphi_\Lambda](s)\big|=O\left(\frac1{|s|}\right),\qquad
|s|\to+\infty,\,\,\,\Re s>-2,
\end{equation*}
so that according to (\ref{4-4.7}), we have 
\begin{equation*}
\big|\lap[\varphi-\varphi_\Lambda](s)\big|=O\left(\frac1{\sqrt{|s|}}
\right),\qquad\text{as}\,\,\, |s|\to+\infty,\,\,\,\Re s>-2.
\end{equation*}
Using some basic complex analysis, we get the intermediate growth control
\begin{equation*}
\big|\lap[\varphi-\varphi_\Lambda](s)\big|
=O\left(\frac1{|s|^{(3+\Re s)/2}}
\right),\qquad\text{as}\,\,\, |s|\to+\infty,\,\,\,-2<\Re s<-1.
\end{equation*}
In particular, by the Plancherel identity, the function
\begin{equation}
t\mapsto e^{(2-\epsilon)t}
\big(\varphi_{\varrho,\theta}(t)-\varphi_\Lambda(t)\big)
\label{approx-1}
\end{equation}
is in $L^2(\R_+)$ for each positive $\epsilon$. Together with 
\begin{equation}
\varphi_{\varrho,\theta}(t)-\varphi_\Lambda(t)=O\big(t^{-3/2}e^{-t}\big)
\qquad\text{as}\,\,\,t\to+\infty,
\label{approx-1-1}
\end{equation}
this gives us rather good asymptotic information regarding the behavior of 
the difference function $\varphi_{\varrho,\theta}-\varphi_\Lambda$.

\subsection{The asymptotic formula}

\noindent{\bf Generalized Taylor series for $\Lambda$.} 
The function $\Lambda(x)$ has a convergent series expansion
$$\Lambda(x)=\sqrt{x-1}\Big(\lambda_0+\lambda_1(x-1)+\lambda_2(x-1)^2+\ldots
\Big),\qquad 1\le x<2.$$
It is well-known in asymptotic analysis (see \cite{BH}) 
that this leads to an asymptotic 
formula for its Fourier-Laplace transform $\varphi_\Lambda$,
$$\varphi_\Lambda(t)=\sum_{j=0}^{N}\lambda_j\,\Gamma\left(j+\frac32\right)\,
t^{-j-3/2}e^{-t}+O\big(t^{-N-5/2}e^{-t}\big)$$
as $t\to+\infty$. In view of (\ref{approx-1}) and (\ref{approx-1-1}),
we get the same type of expansion for $\varphi_{\varrho,\theta}$:
\begin{equation}
\varphi_{\varrho,\theta}(t)=
\sum_{j=0}^{N}\lambda_j\,\Gamma\left(j+\frac32\right)\,
t^{-j-3/2}e^{-t}+O_{L^p}\big(t^{-N-2}e^{-t}\big),
\label{eq-asympexp}
\end{equation}
where by $O_{L^p}(1)$ we mean an expression with bounded norm in $L^p(\R_+)$,
and $2<p<+\infty$ is arbitrary. 
The coefficients $\lambda_j$ may be read off from the identity 
(\ref{deflambda}); for instance, if we use in addition the identity
(\ref{4-4.4.1}), we obtain 
\begin{multline}
\lambda_0=\frac{\theta\mu\,e^{-\mu\varrho}}{\pi(1-\varrho)}
-\frac1{\pi}\int_0^{+\infty}\big|\varphi'_{\varrho,\theta}(t)\big|\,
\varphi_{\varrho,\theta}(t)\,e^{-\varphi(t)\varrho}\,e^{t}\,\diff t\\
=\frac1{\pi}\int_0^{+\infty}\big|\varphi'_{\varrho,\theta}(t)\big|\,
\big(\mu-\varphi_{\varrho,\theta}(t)\big)\,
e^{-\varphi(t)\varrho}\,e^{t}\,\diff t.
\label{eq-lambda0}
\end{multline}
By (\ref{eq-asympexp}) with $N=1$, this leads to
$$\varphi_{\varrho,\theta}(t)=
\beta_1\,t^{-3/2}e^{-t}+O_{L^p}\big(t^{-2}e^{-t}\big),$$
for any $p$, $2<p<+\infty$, with 
$\beta_1=\Gamma(3/2)\,\lambda_0$. This is the formula alluded to in the
introduction, subsection \ref{subsec-1.3}.

\section{Further topics}

\noindent{\bf Maximal meromorphic extensions.} 
It is of interest to apply the general scheme (\ref{3-8.3}) not just
to $N=1$, but also $N=2,3,4,\ldots$. Unfortunately, this approach seems
to run into difficulty after a few steps. In any case, it is possible to
show that $\lap[\varphi_{\varrho,\theta}](s)$ has an analytic extension
to the half-plane $\Re s>-4$ minus the slit $]-4,-1]$ along the real axis. 
One possible approach to get an explicit expression for 
$\lap[\varphi_{\varrho,\theta}](s)$ would be to obtain the maximal meromorphic
continuation to a Riemann surface sheeted over the complex plane. If the
structure of the Riemann surface would happen to be simple, it might be 
possible to express $\lap[\varphi_{\varrho,\theta}](s)$ in terms of the 
poles of the maximal extension. This would then also give an explicit way
to compute the free boundary function $\varphi_{\varrho,\theta}(t)$ itself.
\medskip

\noindent{\bf Small values of $\theta$.} 
The parameter $\theta$ was introduced mainly for the purpose of having
the opportunity to choose $\theta$ close to $0$, since everything is 
completely understood for $\theta=0$. It might then, in a second step, be 
possible to extend the analysis to general values of $\theta$.

As we differentiate the balayage equation (\ref{eq-BE-3}) with respect to
the parameter $\theta$, we have
\begin{multline*}
\int_0^{+\infty}\frac{\diff\varphi_{\varrho,\theta}}{\diff\theta}(t)\,
e^{-\varphi_{\varrho,\theta}(t)[\varrho+\sqrt{s}]}
e^{-(s-1)t}\,\diff t\\
=\left[\frac{\diff\mu}{\diff\theta}-\frac1{1-\varrho}
+(\varrho+\sqrt{s})\,\frac{\theta}{1-\varrho}\,
\frac{\diff\mu}{\diff\theta}\right]
\frac{e^{-\mu\,(\varrho+\sqrt{s})}}{s-1},
\end{multline*}
and after another differentiation, we have
\begin{multline*}
\int_0^{+\infty}\frac{d^2\varphi_{\varrho,\theta}}{\diff\theta^2}(t)\,
e^{-\varphi_{\varrho,\theta}(t)[\varrho+\sqrt{s}]}e^{-(s-1)t}\,\diff t\\
-(\varrho+\sqrt{s})\int_0^{+\infty}
\left(\frac{\diff\varphi}{\diff\theta}(t)\right)^2
\,e^{-\varphi_{\varrho,\theta}(t)[\varrho+\sqrt{s}]}e^{-(s-1)t}\,\diff t\\
=\Bigg\{\frac{d^2\mu}{\diff\theta^2}+(\varrho+\sqrt{s})\bigg[\frac{\theta}
{1-\varrho}\,\frac{d^2\mu}{\diff\theta^2}+\frac2{1-\varrho}\,
\frac{\diff\mu}{\diff\theta}-
\left(\frac{\diff\mu}{\diff\theta}\right)^2\bigg]\\
-(\varrho+\sqrt{s})^2
\frac{\theta}{1-\varrho}\,\frac{\diff\mu}{\diff\theta}
\Bigg\}\frac{e^{-\mu\,(\varrho+\sqrt{s})}}{s-1}.
\end{multline*}
For $\theta=0$, $\varphi_{\varrho,\theta}(t)\equiv0$, and as we apply the 
first identity,
we realize that we also have
$$\frac{\diff\varphi_{\varrho,\theta}}{\diff\theta}(t)\bigg|_{\theta=0}=0,
\qquad 0<t<+\infty.$$
Inserting $\theta=0$ also into the second, we get
$$\int_0^{+\infty}\frac{d^2\varphi}{\diff\theta^2}(t)\bigg|_{\theta=0}
e^{-(s-1)t}\,\diff t=\frac1{(1-\varrho)^2}\,\frac1{\sqrt{s}+1},
\qquad \Re s>1,$$
and, as a consequence,
$$\frac{d^2\varphi_{\varrho,\theta}}{\diff\theta^2}(t)\bigg|_{\theta=0}=
\frac{\pi^{-1/2}}{2\,(1-\varrho)^2}
\int_t^{+\infty}\tau^{-3/2}e^{-\tau}\diff\tau,\qquad 0<t<+\infty.$$
This gives us some understanding of the behavior of 
$\varphi_{\varrho,\theta}(t)$ for small $\theta$ and fixed $\varrho$ and $t$.

\bigskip

\noindent Department of Mathematics, Royal Institute of Technology, 
S--100 44 Stockholm, Sweden;

\noindent E-mail: {\tt haakanh@math.kth.se}.

\enddocument